\documentclass[12pt]{article}
\usepackage{amssymb,latexsym, amsmath, amscd, array, graphicx}
\makeatletter
\date{}

\newtheorem{theorem}{Theorem}[section]

\newtheorem{proposition}[theorem]{Proposition}

\newtheorem{corollary}[theorem]{Corollary}

\newtheorem{problem}[theorem]{Problem}

\newcommand{\rp}{{\Bbb R}P}
\newcommand{\edim}{{\rm e}$-${\dim}}

\newcommand{\z}{{\Bbb Z}}
\newcommand{\I}{{\Bbb I}}

\newcommand{\bb}{{\Bbb B}}
\newcommand{\re}{{\Bbb R}}

\newcommand{\s}{{\Bbb S}}

\newcommand{\lo}{\longrightarrow}

\newcommand{\black}{{\blacksquare}}

\begin{document}

\title{Extensions of maps to $M(\z_m,1)$}

\author{Jerzy Dydak and  Michael Levin}

\maketitle
\begin{abstract} 
We show that a Moore space $M(\z_m,1)$ is an absolute extensor for
finite dimensional metrizable spaces of cohomological dimension $\dim_{\z_m} \leq 1$.
\bigskip
\\
{\bf Math. Subj. Class.:}  55M10 (54F45, 55N45)
\end{abstract}
\begin{section}{Introduction}
\label{introduction}
All  spaces are assumed to be   metrizable. A map means
a continuous function and a compactum means a compact metrizable space. 
By   cohomology we always mean the \v Cech
cohomology. Let $G$ be an abelian group. The  cohomological dimension
 $\dim_GX$
of a space $X$ with respect to the coefficient group $G$ does not exceed $n$, $\dim _G X \leq n$ if $H^{n+1}(X,A;G)=0$ for every closed $A
\subset X$.  Note that this condition implies that
$H^{n+k}(X,A;G)=0$ for all $k\ge 1$ \cite{Ku},\cite{DrArxiv}.
Thus, $\dim _G X =$ the smallest integer $n\geq 0$ satisfying
$\dim _G X \leq n$ (provided it exists), and $\dim _G X = \infty $ if such an integer
does not exist.

Cohomological dimension is characterized by the following basic property:
$\dim_G X \leq n$ if and only  for every closed
$A \subset X$ and a map $f : A \lo K(G,n)$,
$f$ continuously extends over $X$ where $K(G,n)$ is the Eilenberg-MacLane complex
of type $(G,n)$
(we assume that $K(G,0)=G$ with discrete topology and $K(G, \infty)$ is a singleton).
This extension characterization of Cohomological
Dimension gives a rise to  Extension  Theory (more general than
Cohomological Dimension Theory)
 and  the notion of Extension Dimension.  
 The {\em extension dimension} of a space $X$ is said
to be dominated by a CW-complex $K$, written $\edim X \leq K$, if
every map $f : A \lo K$ from a closed subset $A$ of $X$
continuously extends over $X$. Thus $\dim_G X \leq n$ is equivalent
to $\edim X \leq K(G,n)$ and $\dim X \leq n$ is equivalent to 
$\edim X \leq S^n$. The property $\edim X \leq K$ is also denoted
by $X\tau K$ and it is also referred to as $K$ being an absolute extensor of $X$.

The following theorem shows a  close connection between extension
and cohomological dimensions.

\begin{theorem}[Dranishnikov Extension Theorem]
\label{dranishnikov-extension-theorem}
Let $K$ be a  CW-complex and
  $X$ a  metrizable space. Denote by $H_*(K)$ the reduced integral homology
  of $K$. Then

(i)
$\dim_{H_n(K)} X \leq n$ for every $n\geq 0$
 if
$\edim X \leq K$;

(ii) $\edim X \leq K$ if $K$ is simply connected, $X$ is finite dimensional
and
$ \dim_{H_n(K)} X \leq n$ for every $n\geq 0$.
\end{theorem}
Theorem \ref{dranishnikov-extension-theorem} was proved in \cite{DrExt} for the compact case and extended in \cite{Dy} to the metrizable case.

Let $G$ be an abelian group.  
We always assume that 
a Moore space $M(G,n)$ of  type $(G,n)$ is an $(n-1)$-connected CW-complex whose reduced integral homology is concentrated
in dimension $n$ and equals $G$.
Theorem \ref{dranishnikov-extension-theorem} implies
that for a finite dimensional metrizable space $X$ and $n>1$,  $\dim_G X \leq n$
if and only if $\edim X \leq M(G,n)$. The main open problem for $n=1$ is:
\begin{problem}
\label{problem}
Let $G$ be an abelian group and let $M(G,1)$ be a Moore space whose
fundamental group is abelian. Is $M(G,1)$ an absolute extensor for 
finite dimensional metrizable spaces  of $\dim_G \leq 1$?
\end{problem}

This problem was affirmatively answered  in \cite{dydak-levin-2}
for $M(\z_2,1)=\rp^2$. In this paper we extend the result of \cite{dydak-levin-2}
to
Moore spaces $M(\z_m, 1)$. In this particular case we choose a specific model and 
by $M(\z_m,1)$ we  mean
the space obtained by attaching a  disk  to a circle by an $m$-fold covering map
of the disk boundary. Our main result is:
\begin{theorem}
\label{main-theorem}
The Moore space $M(\z_m,1)$ is an absolute extensor for
finite dimensional metrizable spaces of cohomological dimension mod $m$ at most $1$.
\end{theorem}

The case of metrizable spaces of $\dim \leq 3$ in Theorem \ref{main-theorem}
was independently obtained by A. Nag\'orko by generalizing the approach
of \cite{first-projective} to $3$-dimensional lens spaces.

\end{section}
\begin{section}{Preliminaries}
\label{preliminaries} 
In this section we present a few general notations and facts  
that will be used later.

For a CW-complex $L$ we denote by $L^{[k]}$ the $k$-skeleton
of $L$.

Let  $A$ and $B$ be compact spaces and 
 $A' \subset A$ are  $B'\subset B$ closed subsets.
We denote by $\frac{A \times B}{A' \times B'}$ the quotient
space of  $A\times B$  
by the partition consisting of   the singletons
of $(A\times B) \setminus (A'\times  B') $ and the sets 
$\{ a \} \times B', a \in A'$. Clearly, 
the spaces $\frac{A' \times B}{A' \times B'}$ and
$\frac{A \times B'}{A' \times B'}$
  can be considered  as
subspaces of the space  $\frac{A\times B}{A'\times B'}$.
In a similar way we define for closed subsets
 $A'\subset A'' \subset A$ and
$B'\subset B'' \subset B$ the space 
$$ \frac{A \times B}{ A''\times B' \cup A ' \times B''}$$
 as the quotient space of $A\times B$ by  the partition 
 consisting of 
 the sets $\{a \} \times B'$ for $ a \in A''\setminus A'$,
 the sets $\{ a \} \times B''$ for $a \in A'$ and the singletons 
 not contained in the sets listed before.

\begin{proposition}
\label{composition}
Let $f,g : M(\z_m,1) \lo M(\z_m,1)$ be maps 
inducing the zero-homomorphism of the fundamental group.
Then $f\circ g$ is null-homotopic.
\end{proposition}
{\bf Proof.} Note that the universal cover ${\tilde M}$ of $M(\z_m,1)$
is homotopy equivalent to a bouquet of $2$-spheres and the maps
$f$ and $g$ lift to ${\tilde M}$. Thus the map  $f \circ g$
factors 
through maps ${\tilde M} \lo M(\z_m,1) \lo {\tilde M}$ whose
composition
induces the zero-homomorphism
of $H_2({\tilde M})$  since  $H_2(M(\z_m, 1))=0$. Hence the composition
${\tilde M} \lo M(\z_m,1) \lo {\tilde M}$ is a null-homotopic map resulting in $f\circ g$ being null-homotopic.
$\black$.

\begin{proposition}
\label{trivial-bundle}
Let $T$ be a  $G$-bundle
over $M(\z_m,1)$ and $M_0$ a singleton in  $ M(\z_m,1)$.  If the structure group $G$ of the bundle
is arcwise connected 
 then $T$ is trivial over
$M(\z_m,1) \setminus M_0$.
\end{proposition}
{\bf Proof.} Take a sufficiently fine triangulation 
of $M(\z_m,1)$ 
and observe that there is a retraction 
$r : M(\z_m,1) \setminus M_0 \lo L$ 
to a $1$-dimensional subcomplex $L$ of $M(\z_m,1)$ such that
$r$ can decomposes into the composition of retractions 
that move the points of $M(\z_m,1)$ only inside  small
sets over which $T$ is trivial. Then $r$ induces a bundle
map from $T$ over $M(\z_m,1) \setminus M_0$
to $T$ over $L$. Note that every $G$-bundle
over a one-dimensional simplicial complex is trivial
if $G$ is arcwise connected. Thus $T$ over $L$ is trivial 
and hence $T$ over $M(\z_m,1) \setminus M_0$
is trivial as well.
$\black$
\\\\
The following  two propositions are simple exercises left to the reader.
\begin{proposition}
\label{embedding}
Let $T$ be a ball bundle over a metrizable space $L$,
$T_0$  the fiber of $T$ over a point in $ L$ and
$U$  
a neighborhood of $T_0$ in $T$. Then 
$T/T_0$ embeds into $T$ so that $T\setminus (T/T_0)\subset U$
and the projection of $T/T_0$ 
to $L$ coincides with the projection of $T$ to $L$ restricted
to $T/T_0$.
\end{proposition}

\begin{proposition}
\label{lifting}
 Let $X$ be a metrizable space, 
$g : K \lo L$ a map of a CW-complex $K$ to a simplicial complex
$L$ (with the CW topology) such that for every simplex $\Delta$ of $L$ we have that
$g^{-1}(\Delta)$ is a subcomplex of $K$ and
$\edim X \leq g^{-1}(\Delta)$. Then 
a map $f : F \lo K$ from a closed subset $F$
of $X$  extends over $X$ if $g\circ f : F \lo L$
extends over $X$.
\end{proposition}
We will also need
\begin{proposition} 
{\rm (\cite{dydak-levin-2})}
\label{suspension}
Let $K$, $L$ and $M$ be finite CW-complexes, $L_0$ a singleton in $L$,
$X$ a metrizable space and  $F$ closed subset of $X$ 
 such that $L$ is connected,
$L$ admits a simplicial structure and $\edim X \leq \Sigma K$.

(i) If a map  $f : F \lo \frac{L\times K}{L_0 \times K}$  followed
by the projection of $\frac{L\times K}{L_0 \times K}$ to $L$
extends over $X$ then $f$ extends over $X$ as well.

(ii) If $f : F \lo L\times K$ and $g : L\times K \lo M$ are maps
such that $f$ followed by the projection of $L \times K$ to $L$
extends over $X$
and $g$ is null-homotopic on $L_0 \times K$ then $f$ followed by $g$ 
extends over $X$ as well.

\end{proposition}
For the reader's convenience let us outline the proof
of Proposition \ref{suspension}. Fix a triangulation
of $L$ for which $L_0$ is a vertex. Observe that
the projection of $\frac{L\times K}{L_0 \times K}$ to $L$
factors up to homotopy through the space 
$\frac{L\times K}{L^{[0]}\times K}$ 
where
$L^{[0]}$ is the $0$-skeleton of $L$ with respect
to the triangulation of $L$. Also observe that
for every $n$-simplex $\Delta$ the space 
$\frac{\Delta \times K}{\Delta^{[0]} \times K}$ is homotopy
equivalent to the wedge of $n$ copies of $\Sigma K$
and hence $\edim X \leq \frac{\Delta \times K}{\Delta^{[0]} \times K}$.
Then, by Proposition \ref{lifting}, we have that $f$ in (i) followed
by the projection of $L \times K$ to 
$\frac{L \times K}{ L^{[0]}\times K}$ extends over $X$.
Hence $f$ extends over $X$ as well and (i) is proved.
Note that in (ii) the map $g \times f$ factors up to homotopy
through
$ \frac{L \times K}{ L_0\times K}$ 
and hence (ii) follows from (i).

\end{section}
\begin{section}{Lens spaces}
\label{lens-spaces}

By $\re^n$, $\bb^n$, $\s^n$ we denote the $n$-dimensional Euclidean 
space, the unit ball in $\re^n$, and the unit sphere in $\re^{n+1}$
respectively. A topological $n$-sphere is denoted by $S^n$ with
$S^0$ being a singleton.
We usually  assume that that $\re^m \subset \re^k$ if $m \leq k$.
Thus we will use the subscript $\perp$
to write $\re^n=\re^m \oplus \re^k_\perp$ for $n=m+k$ in order 
to emphasize that $\re^m$ and  $\re^k_\perp$ are  not subspaces
of each other.

Recall that for  a CW-complex $L$ we denote by $L^{[k]}$ the
$k$-skeleton  of $L$.   For a covering space   $\tilde L$ of $L$
we will consider $\tilde L$ with the CW-structure induced 
by the CW-structure of $L$ and hence we have that the $k$-skeleton 
${\tilde L}^{[k]}$ of $\tilde L$ is the preimage  of $L^{[k]}$ under the covering map.

In the proof of Theorem \ref{main-theorem}
we will use the infinite dimensional lens  space $L_m$ as a  model of $K(\z_m,1)$.
 Let
us remind the construction of $L_m$.
 Decompose  $\re^{2n}$ into the direct sum
of $n$ two-dimensional coordinate planes $\re^2$ 
and consider the orthogonal transformation $\theta$
of  $\re^{2n}$ induced by rotating
counterclockwise
each $\re^2$ in the decomposition of $\re^{2n}$ by the angle $2\pi /m$.
Thus $\z_m=\{\theta,\theta^2, \dots, \theta^m \}$ acts 
on $\re^{2n}$ by orientation preserving 
orthogonal transformations and 
 $\z_m$ acts freely on the unit sphere $\s^{2n-1}$ of $\re^{2n}$.
We will refer to $\theta$ as the generating transformation of $\z_m$.
Denote $L_m^{[2n-1]}=\s^{2n-1}/\z_m$. Representing $\re^{2n+2}$ as  
 $\re^{2n+2}=\re^{2n} \oplus \re^2_\perp$ we can regard $\s^{2n-1}$ as a subset of
$\s^{2n+1}$ and $L_m^{[2n-1]}$ as a subset of 
$L_m^{[2n+1]}$. 
The infinite dimensional lens space $L_m$ is 
defined as $L_m=$dirlim $L_m^{[2n-1]}$.
Clearly  $\z_m$ freely acts on 
  $\s^\infty$=dirlim $\s^{2n-1}$  and
 $L_m=\s^\infty /\z_m$. Thus we have that
$L_m =K(\z_m,1)$  since $\s^\infty$ is contractible.  

The CW-structure of $L_m$ is defined  so that $L_m$ has
only one cell in  each dimension, see \cite{hatcher}.
The CW-structure of $L_m$ agrees with 
our previous notation since 
$L^{[2n-1]}_m$  is  indeed the $(2n-1)$-skeleton
of $L_m$.   Set  ${\tilde L}_m =\s^\infty$ and  ${\tilde L}^{[2n-1]}_m=\s^{2n-1}$.
The preimage  ${\tilde L}^{[2n]}_m $  of ${ L}^{[2n]}_m $ under the projection
of $\s^{2n+1}$ to ${ L}^{[2n+1]}_m $
can be described as follows.
Represent $\re^{2n+2}$ as the direct sum $\re^{2n+2}=\re^{2n} \oplus \re^{2}_\perp$
of orthogonal coordinate subspaces invariant under the action of $\z_m$
 on $\re^{2n+2}$ and let 
  $\s^{2n-1}$  and $\s^{1}_\perp$  be the unit  sphere and the unit circle   in 
$\re^{2n}$ and $\re^{2}_\perp$ respectively. Take a point  $ a \in  \s^1_\perp$
and  consider the $(2n+1)$-dimensional linear subspace of $\re^{2n+2}$
containing $\re^{2n}$ and the point $a$, and in this subspace consider the unit
sphere $S^{2n}_a$. Then $\s^{2n-1}$ divides $S^{2n}_a$ into two hemispheres
and denote by $C_a$ the hemisphere containing the point $a$. It is clear
that $g C_a=C_{ga}$ for $g \in \z_m$. 
Fix any orbit $A$ of $\z_m$ in $\s^1_\perp$.
The space ${\tilde L}^{[2n]}_m$ is
defined as the union of 
${\tilde L}^{[2n-1]}_m=\s^{2n-1}$ with the $(2n)$-dimensional hemispheres
 $C_a, a\in A$ which
are defined to be the $(2n)$-cells of ${\tilde L}^{[2n]}_m$. Clearly 
${\tilde L}^{[2n]}_m$
is invariant under the action of $\z_m$ on $\re^{2n+2}$. Then
$L^{[2n]}_m$ is defined as the orbit space
$L^{[2n]}_m={\tilde L}^{[2n]}_m/\z_m$ and it is obvious that
$L^{[2n]}_m$ is obtained from $L^{[2n-1]}_m$ by attaching one $(2n)$-cell.

We will call the models of $L^{[2n-1]}_m$  and $L^{[2n]}_m$ described above 
the covering models. Note that ${\tilde L}^{[k]}_m$ is the universal cover
of $L^{[k]}_m$ for $k>1$ and  ${\tilde L}^{[1]}_m=\s^1$ is an $m$-fold cover
of $L^{[1]}_m=S^1$. The $0$-skeleton $L^{[0]}_m$ of $L_m$ 
is a singleton
in $L^{[1]}_m$. Note that $L^{[2]}_m=M(\z_m,1)$.

The space $L^{[2n]}_m$ can be also described in the following way.
Consider the unit ball $\bb^{2n}$ in $\re^{2n}$. Then $L^{[2n]}_m$
is the quotient space of $\bb^{2n}$ under the action of $\z_m$ on 
 $\partial \bb^{2n}=\s^{2n-1}$. By this we mean the  quotient space
whose equivalence classes are the orbits of the points in $\s^{2n-1}$
and the singletons in $\bb^{2n} \setminus \s^{2n-1}$. We will refer
to such representation of $L^{[2n]}_m$ as the ball model 
of $L^{[2n]}_m$.

The space $L^{[2n+1]}_m$  also admits  a similar  description.
Represent  $\re^{2n+1}=\re^{2n} \oplus \re_\perp$, consider
the unit sphere $\s^{2n-1}$ in $\re^{2n}$ and the action of $\z_m$
on $\re^{2n}$.
Consider the unit ball $\bb^{2n+1}$ in $\re^{2n+1}$ and 
define an equivalence relation on  $\bb^{2n+1}$ with the equivalence
classes   to be: the orbits of the action of  $\z_m$
on $\s^{2n-1}$, the singletons of $\bb^{2n+1} \setminus \partial \bb^{2n+1}$
and the sets $\{ (x,t), (\theta x, -t) \}$ where $\theta$
 is the generating transformation of $\z_m$ 
and
$(x,t)\in \re^{2n+1}=\re^{2n} \oplus  \re_\perp$ such that
$(x,t) \in \partial \bb^{2n+1}$ and $t <0$. Then $L^{[2n+1]}_m$ is 
the quotient space of $\bb^{2n+1}$ under this equivalence relation.
Similarly we  refer to such representation of $L^{[2n+1]}_m$ as 
the ball model of $L^{[2n+1]}_m$.

\end{section}
\begin{section}{Extensions of maps to Lens spaces}
\label{extensions-to-lens-spaces}
In this section we prove two auxiliary propositions. 
By a Moore space $M(\z_m, k)$ we mean a space 
obtained by attaching a $(k+1)$-ball to a $k$-sphere  $S^k$
by a map degree $m$ from
the ball boundary to $S^k$ 
and we denote the $k$-sphere $S^k$ in $M(\z_m,k)$ by $\partial M(\z_m, k)$.

\begin{proposition}
\label{extension-odd}
Let $\psi: S^1 \times S^{2n-1} \lo L^{[2n]}_m, n\geq 1$ be a map such 
that $\psi$ restricted to  $S^1 \times S^0$ generates
 the fundamental group of $L^{[2n]}_m$ and
 $\psi$ restricted to 
$S^0 \times S^{2n-1}$ is null-homotopic.
Then $\psi$ considered as a map from
 $ S^1 \times \partial M(\z_m, 2n-1)$  
 extends over $S^1 \times M(\z_m, 2n-1)$.
\end{proposition}
{\bf Proof.} Replacing $\psi$ by a homotopic map
assume that $\psi$ factors through 
$\frac{S^1 \times S^{2n-1}}{S^0 \times S^{2n-1}}$.
 Represent 
 $S^1$ as the quotient space of $I=[0,1]$ under the projection 
sending the end-points $\partial I$ of $I$  to $S^0$ and consider
the induced  projection from  the $(2n)$-sphere
$\frac{I \times S^{2n-1}}{\partial I \times S^{2n-1}}$
to $\frac{S^1 \times S^{2n-1}}{S^0 \times S^{2n-1}}$.
Then this projection followed by the map induced by 
$\psi$ from $\frac{S^1 \times S^{2n-1}}{S^0 \times S^{2n-1}}$
to $L^{[2n]}_m$ 
 lifts to a map
 ${ \psi}_I : \frac{I \times S^{2n-1}}{\partial I \times S^{2n-1}}
\lo {\tilde L}^{[2n]}_m$ to the universal cover ${\tilde L}^{[2n]}_m$
of $ L^{[2n]}_m$. Denote by $g \in \z_m$  the element of  the fundamental
group $\z_m$ of $L^{[2n]}_m$  represented by the map $\psi$ restricted to
$S^1 \times S^0$ with $S^0 \times S^0$ and  $\psi(S^0\times S^0)$ being
the base points in  $S^1 \times S^0$  and $L^{[2n]}_m$ respectively and
recall that $g$ is a generator of $\z_m$.

Represent $M(\z_m,2n-1)$ as the quotient space
of a $(2n)$-ball $B$ under the projection from  $B$ to $M(\z_m, 2n-1)$
sending $\partial B=S^{2n-1}$ to $\partial M(\z_m,2n-1)=S^{2n-1}$ 
by a map of degree  $m$. Consider the induced projection
from $\frac{I \times B }{\partial I\times  B}$ 
to $\frac{ I \times M(\z_m,2n-1)}{\partial I \times M(\z_m,2n-1)}$ and
 denote by 
$\psi_B: \frac{I \times \partial B }{\partial I \times \partial B} \lo 
{\tilde L}^{[2n]}_m$ this projection restricted
to $\frac{I \times \partial B }{\partial I \times \partial B}$ and 
 followed by the map ${ \psi}_I$.
Then 
the problem of extending $\psi$ reduces to the problem of
extending $\psi_B$  to a map
$\psi'_B : \frac{I \times  B }{\partial I  \times \partial B}\lo 
{\tilde L}^{[2n]}_m$
so that for every $x \in B /\partial B$  and
 $(0,x),(1,x) \in \partial I \times (B/\partial B)=
\frac{\partial I \times B}{\partial I \times \partial B}
\subset \frac{I \times  B }{\partial I  \times \partial B}$
we have that $\psi'_B (1,x) =g(\psi'_B(0,x))$
with the element $g$ of the fundamental group of $L^{[2n]}_m$
being considered as
acting on ${\tilde L}^{[2n]}_m$. 

Note that 
$S^{2n}_\#=
\frac{I \times \partial B }{\partial I \times \partial B}$ is a $(2n)$-sphere
and 
$\psi_\# =\psi_B | S^{2n}_\# : S^{2n}_\# \lo {\tilde L}^{[2n]}_m$
 factors through  a map of degree $m$ from
$S^{2n}_\#$ to
$S^{2n}_M=
\frac{ I \times \partial M(\z_m,2n-1)}{\partial I \times \partial M(\z_m,2n-1)}$.
Also note that under 
the projection of $I \times B$ to 
$ \frac{I \times  B }{\partial I  \times \partial B}$ 
the $(2n)$-sphere  $S^{2n}_+=\partial (I \times B)$
 goes to the space
$ \frac{I \times  B }{\partial I  \times \partial B}$ which is
the union of the spheres $S^{2n}_0 =\{0 \} \times (B / \partial B)$,
$S^{2n}_1 =\{1 \} \times (B / \partial B)$ and
$S^{2n}_\#$
so that $S^{2n}_0$ and $S^{2n}_1$ are disjoint and each of them intersects
$S^{2n}_\#$ at only one point. And finally note   that
${\tilde L}^{[2n]}_m$ is homotopy equivalent to a bouquet  of $(2n)$-spheres.
Then the problem of extending $\psi_B$ to $\psi'_B$ boils down
to constructing a map $\psi_0 : S^{2n}_0 \lo {\tilde L}^{[2n]}_m$
so that $(\psi_0)_*(\alpha)-g_*((\psi_0)_*(\alpha))+(\psi_\#)_*(\beta)=0$
in the homology group $H_{2n}({\tilde L}^{[2n]}_m)$ with
$\alpha$ and $\beta$ being the generators  of $H_{2n}(S^{2n}_0)$
and $H_{2n}(S^{2n}_\#)$ determined by the orientations of
$S^{2n}_0$ and $S^{2n}_\#$ induced by an orientation
of $S^{2n}_+$. 

Recall that $\gamma=(\psi_\#)_*(\beta)$ is divisible by $m$
since $\psi_\#$ factors through a map of degree $m$ of
a $(2n)$-sphere
and $g$ comes from an orientation preserving
orthogonal transformation  of the $(2n+1)$-sphere 
${\tilde L}^{[2n+1]}_m$. Consider the cellular homology of 
${\tilde L}^{[2n]}_m$, fix an oriented $(2n)$-cell $C$ of 
${\tilde L}^{[2n]}_m$ and index the  $(2n)$-cells
$C_1,C_2, \dots, C_m$ of ${\tilde L}^{[2n]}_m$  so that $C_i=g^i(C), 1\leq i \leq m$.
Let $ \gamma_1 C_1 +\gamma_2 C_2 +\dots +\gamma_m C_m,
\gamma_1+ \dots +\gamma_m=0,
 \gamma_i \in \z$, be the cycle representing $\gamma$
 and
 $y_1 C_1 +\dots+ y_m C_m, y_1 + \dots + y_m=0, y_i \in \z,$
 the cycle representing $y=(\psi_0)_*(\alpha)$. 
 Then $g_*(y) $ is represented
 by the cycle $y_m C_1+y_1 C_2 +\dots + y_{m-1} C_m$ and
  we arrive
 at the system of linear equations over $\z$: 

$$\begin{cases}
\gamma_1+ \dots +\gamma_m=0\\
y_1 + \dots +y_m=0\\
y_1-y_m+\gamma_1=0\\
y_2-y_1+\gamma_2=0\\
\dots \\
y_m -y_{m-1}+\gamma_m=0.
\end{cases}$$
Representing  $y_{m}=-y_1-y_2-\dots -y_{m-1}$    get
$$\begin{cases}
2y_1 + y_2 +\dots + y_{m-1}+\gamma_1=0\\
y_2-y_1+\gamma_2=0\\
\dots \\
y_{m-1} -y_{m-2}+\gamma_{m-1}=0.
\end{cases}$$
Eliminating $y_1, \dots, y_{m-2}$ from the first equation get 
$$y_{m-1}= -\frac{1}{m}(\gamma_1+2\gamma_2 +3\gamma_3 + \dots + (m-1)\gamma_{m-1})$$
 and  find $y_{m-2}, y_{m-3}, \dots, y_1$  from the remaining equations.
Recall that $\gamma$ is divisible by $m$ and hence 
$\gamma_1, \dots,\gamma_m$ are divisible by $m$ as well. 
Thus we conclude that the system is solvable over $\z$.
Set $\psi_0$ to be a map with  $(\psi_0)_*(\alpha)=y$ and 
the proposition
follows. $\black$

\begin{proposition}
\label{extension-even}
Let 
$\psi : 
\frac{L^{[2]}_m\times S^{2n-1}}{ L^{[0]}_m\times S^{2n-1}} 
\lo L^{[2n+1]}_m, n\geq 1,$ be a map. Then
$\psi$ considered as 
a map from 
$\frac{L^{[2]}_m\times \partial M(\z_m,2n-1)}
{L^{[0]}_m \times \partial M(\z_m,2n-1)}$
extends over
$\frac{L^{[2]}_m\times  M(\z_m,2n-1)}
{L^{[0]}_m \times  M(\z_m,2n-1)}$.

\end{proposition}
{\bf Proof.} 
Note that  $L^{[1]}_m \subset L^{[2]}_m=
\frac{ L^{[2]}_m \times S^0}{L^{[0]}_m \times S^0}
\subset 
\frac{ L^{[2]}_m \times S^{2n-1}}{L^{[0]}_m \times S^{2n-1}}$,
 denote by $g \in \z_m$ the element of the fundamental group
of $L^{[2n+1]}_m$ represented by $\psi$ restricted to 
the circle $L^{[1]}_m$ and consider 
 $g$ as an orthogonal transformation
 acting on  the universal cover
 $S^{2n+1}= {\tilde L}^{[2n+1]}_m$ of $L^{[2n+1]}_m$.
 
 Let the projection $p_I : [0,1] \lo L^{[1]}_m$
 send the  end points $\partial I$ of $I$  to $L^{[0]}_m$.
 This projection induces a projection of 
 the $(2n)$-sphere
$\frac{I \times S^{2n-1}}{\partial I \times S^{2n-1}}$ to
 $\frac{L^{[1]}_m \times S^{2n-1}}{ L^{[0]}_m\times S^{2n-1}}$
 and this projection followed by $\psi$ lifts
 to a map 
 $\psi_I : \frac{I \times S^{2n-1}}{\partial I \times S^{2n-1}}
 \lo {\tilde L}^{[2n+1]}_m=S^{2n+1}$. Then $\psi_I$
   factors
up to homotopy relative to 
$\frac {\partial I \times S^{2n-1}}{\partial I \times S^{2n-1}}=
\partial I $ 
through the space $\frac { I \times S^{2n-1}}{ I \times S^{2n-1}}
=I$. It implies that $\psi$ factors up to homotopy through
the space $\frac{ L^{[2]}_m \times S^{2n-1}}{L^{[1]}_m \times S^{2n-1}}$.
Thus replacing $\psi$ by a map
from the last space we may assume that 
$\psi : 
\frac{ L^{[2]}_m \times S^{2n-1}}{L^{[1]}_m \times S^{2n-1}} \lo
 L^{[2n+1]}_m$  and look for 
an extension of $\psi$ over the space 
$$\frac{ L^{[2]}_m \times M}
{L^{[0]}_m \times M 
\cup L^{[1]}_m\times \partial M}$$
where we shorten $M(\z_m, 2n-1)$ and $\partial M(\z_m, 2n-1)$
to $M$ and $\partial M$ respectively.

Represent  $L^{[2]}_m$ as the quotient space
of a disk $D$ under the projection  $p_D : D \lo L^{[2]}_m$
which sends $\partial D=S^1$ to $L^{[1]}_m=S^1$ by an $m$-fold map
and denote $D^0=p^{-1}(L^{[0]}_m)$. 

Denote 
$$
K=\frac{D \times M}{D^0\times M \cup 
\partial D \times \partial M} $$
$$\partial_D K =\frac{\partial D \times M}
{D^0 \times M \cup \partial D \times \partial M}\subset K 
$$
$$ \partial_M K=
\frac{D \times S^{2n-1}}{\partial D \times S^{2n-1}}=
\frac{D \times \partial M}{\partial D \times \partial M}\subset K.$$
  The projection from  
$\partial_D K=
\frac{D \times S^{2n-1}}{\partial D \times S^{2n-1}}=S^{2n+1} $ to
$\frac{L^{[2]}_m\times S^{2n-1}}{ L^{[1]}_m\times S^{2n-1}}$ induced
by $p_D$ and followed by $\psi$ lifts to a map
  $\psi_M : \partial_M K 
 \lo 
{\tilde L}^{[2n+1]}_m=S^{2n+1}$.
Consider
a  rotation  of $\partial D$ by the angle $2\pi/m$ under which 
the map $p_D$ restricted  to $\partial D$ is invariant.
 Then this rotation induces 
a rotation (homeomorphism) $\omega$ of the space $\partial_D K$. 
Thus the problem of extending $\psi$ reduces
to the problem of extending $\psi_M$ to a map
$\psi'_M:  K \lo 
 {\tilde L}^{[2n+1]}_m$ so that
 $\psi'_M (\omega (x))=g(\psi'_M(x))$ for $x \in \partial_D K$.
 
 Note that $\partial_D K$ is the union of  $(2n+1)$-spheres
 $S^{2n+1}_1, \dots, S^{2n+1}_m$  intersecting
 each other at points of $ D^0$ (we consider 
 $D$ as  a natural subset of $K$).
 Also note that $ \partial_M K\cap  \partial_D K=\partial D$
  and 
 $\partial_M K=S^{2n+1}_\#$ is 
  a $(2n+1)$-sphere intersecting
 the spheres $S^{2n+1}_i, 1\leq i \leq m$,
 at points of $\partial D$. Clearly $\partial_D K$ is
  invariant under $\omega$ and 
 and  the spheres $S^{2n+1}_i$ can be indexed 
  so that  $S^{2n+1}_i=\omega^i(S^{2n+1}_m)$.
 
 Consider a projection from a $(2n)$-ball $B$ to $M$
 sending $\partial B$ to $\partial M$ by a map
 of degree $m$. This projection
 induces a projection 
 $p:  D \times B \lo K$
 from the $(2n+2)$-ball  $D \times B$ to $K$ 
 under which  the $(2n+1)$-sphere
 $\partial (D\times B)$ goes  to $\partial_M K \cup \partial_D K$
 so that the sphere $S^{2n+1}_\#$ is covered $m$-times and
 each of the spheres $S^{2n+1}_i,  1\leq i \leq m,$ is covered
 only once. Recall that $S^{2n+1}_i =\omega^i (S^{2n+1}_m)$,
 $\psi_M(\omega(x))=g(\psi_M(x))$ and $g$ is an orientation preserving 
 orthogonal transformation of $S^{2n+1}={\tilde L}^{[2n+1]}_m$.
  Consider the spheres $S^{2n+1}_\#, S^{2n+1}_1, \dots S^{2n+1}_m$
  with the orientation induced
  by an orientation of the sphere $\partial(D \times B)$ and
   define a map $\psi_m : S^{2n+1}_m \lo S^{2n+1}={\tilde L}^{[2n+1]}_m$
   so that $\deg \psi_m=-\deg\psi_M |S^{2n+1}_\#$ and
    $\psi_m$ extends $\psi_M$ restricted
   to $S^{2n+1}_m$. Now define
   $\psi_i=g^i \circ \psi_m \circ \omega^{-i} : 
   S^{2n+1}_i   \lo {\tilde L}^{[2n+1]}_m$. Thus
   we have extended $\psi_M$ over  $\partial_M K \cup \partial_D K$
   so that that the map $p$ restricted to 
   $\partial( D \times B)$ and followed by this extension
   is of degree $0$ and hence extends to a map from $D \times B$ 
   to $ {\tilde L}^{[2n+1]}_m$.
   Clearly the last extension induces a map 
   $\psi'_M : K \lo {\tilde L}^{[2n+1]}_m$ with the required properties
 and the proposition follows. $\black$
 
\end{section}

\begin{section}{Pushing maps off  the $(2n+1)$-skeleton of $L_m$ }
\label{odd-reduction-section}
In this section we will prove
\begin{proposition}
\label{odd-reduction-proposition}
Let $X$ be a metrizable space with $\dim_{\z_m} X \leq 2n-1, n\geq 2$,
and let $f : X \lo  L^{[2n+1]}_m$
be  a map. Then there is a  map 
$f' : X \lo  L^{[2n]}_m$ which coincides with $f$ on $f^{-1}(L^{[2n-1]}_m)$.
\end{proposition}
The proof of Proposition \ref{odd-reduction-proposition} 
is based on a modification
of $L^{[2n+1]}_m$. This modification is defined for $n\geq 1$ and
will be referred to as the basic modification of $L^{[2n+1]}_m$.
We describe 
 this modification in such a  way  and  using such  notations 
 that it can  be used in Section 
\ref{even-reduction-section} 
for constructing 
a similar modification of $L^{[2n+2]}_m$.  

Consider the covering model of $L^{[2n+1]}_m$.
 Let $\re^{2n+2}=\re^{2n}\oplus \re^2_\perp$ and let 
$\s^{2n+1}$, $\s^{2n-1}$ and $\s^1_\perp$
be the unit spheres and the unit circle in $\re^{2n+2}$,
 $\re^{2n}$ and $\re^2_\perp$ respectively.
Fix a sufficiently small $\epsilon >0$ and
take a closed  $\epsilon$-neighborhood ${ {E_\s^1}}$ of $\s^1_\perp$ in $\s^{2n+1}$
such that ${ {E_\s^1}}$ 
does not intersect $\s^{2n-1}$. 
Clearly ${{E_\s^1}}$  is invariant
under the action of $\z_m$ on $\s^{2n+1}$ and ${{E_\s^1}}$ can be considered
as a trivial $(2n)$-ball bundle over $\s^1_\perp$ with respect to 
the group $SO(2n)$ of orientation preserving
orthogonal transformations of a $(2n)$-ball.
The bundle ${E_\s^1}$ over $\s^1_\perp$ can be visualized as follows.
Take a point $a \in \s^1_\perp$ and consider the unit sphere $S^{2n-1}_a$ in the linear
$(2n)$-dimensional subspace of $\re^{2n+2}$ containing $\s^{2n-1}$ and $a$.
Then the closed $\epsilon$-neighborhood of $a$ in  $S^{2n}_a$ will be
the $(2n)$-ball over $a$ in the bundle ${E_\s^1}$. The sphere $\s^{2n-1}$
divides $S^{2n}_a$ into two hemispheres and for the hemisphere $C_a$
containing the point $a$ consider the natural
deformation retraction of $C_a\setminus \{a \}$ to $\s^{2n-1}$ along the shortest arcs 
in $S^{n}_a$
connecting  $a$ with the points of $\s^{2n-1}$. Then this retraction
induces the corresponding deformation retraction
${r_\s^1} : \s^{2n+1} \setminus \s^1_\perp \lo \s^{2n-1}$
which commutes with the transformations of $\z_m$. Note
that from this description of the   bundle 
 it can be   seen that 
the transformations
in $\z_m$ induce bundle maps of ${{E_\s^1}}$.

Let $p : \s^{2n+1} \lo L^{[2n+1]}_m =\s^{2n+1}/\z_m$ be the projection.
Denote  $S^1_\perp =p(\s^1_\perp)=\s^1_\perp /\z_m$ and $E^1=p({E_\s^1})$. 
Then 
$E^1$ is a trivial $(2n)$-ball bundle over $S^1_\perp$  
(since ${E_\s^1}$ is a bundle
with respect to the orientation preserving orthogonal transformations)
and
${r_\s^1}$ induces
the deformation retraction 
$ r^1 : L^{[2n+1]}_m\setminus S^1_\perp\lo L^{[2n-1]}_m$.
Represent $E^1=S^1 \times B$ where $B$ is an $(2n)$-ball and denote
$\partial E^1= S^1\times \partial B=S^1\times S^{2n-1}$. 
By $S^0$ we denote
 a singleton in  a sphere $S^k$. Note that, since  $r^1$ is 
a deformation retraction, $r^1$ sends  
the  circle  $S^1 \times S^0 \subset \partial E^1$ 
to a circle in $L^{[2n-1]}_m$ homotopic to $S^1 \times S^0$
in $L^{[2n+1]}_m$. On the other hand $S^1 \times  S^0$
homotopic to the circle $S^1_\perp$ which represents a generator of the fundamental group
of $L^{[2n+1]}_m$ and hence  represents a generator
of the fundamental group
of $L^{[2n]}_m$ as well. 
Also note that $S^0 \times S^{2n-1} \subset \partial E^1$ is 
contractible in the ball
$E^1 \cap  L^{[2n]}_m$. Thus  $r^1$ restricted
to $S^1 \times S^0 $ and $S^0 \times S^{2n-1}$  and followed by
the inclusion of $L^{[2n-1]}_m$ into $L^{[2n]}_m$
represent a generator  of the fundamental
group of $L^{[2n]}_m$ and 
a null-homotopic map to $L^{[2n]}_m$
 respectively.

By the {\bf basic surgery} of $L^{[2n+1]}_m$  we mean
replacing
 $E=E^1=S^1\times B$ with    $E_M=S^1 \times M(\z_m, 2n-1)$ such that 
$\partial E=S^1 \times \partial B$ is identified with 
 $\partial E_M=S^1\times \partial M(\z_m,2n-1)$ through 
an identification of 
$\partial M(\z_m,2n-1)=S^{2n-1}$ with $\partial B=S^{2n-1}$.
The {\bf basic modification}
$M$ of $L^{[2n+1]}_m$ is the space  obtained from $L^{[2n+1]}_m$  by 
the basic surgery of $L^{[2n+1]}_m$. 
  Clearly $L^{[2n-1]}_m$ remains untouched in  $M$.

The basic surgery of $L^{[2n+1]}_m$ 
 can be even easier  described  in the ball
model of $L^{[2n-1]}_m$. 
In this model   the set ${E^1}$ is represented by 
 the closed $\epsilon$-neighborhood of $ \re_\perp \cap \bb^{2n+1} $ in $\bb^{2n+1}$ and
the retraction $r^1$ is represented by the natural  retraction from 
$\bb^{2n+1}\setminus \re$ to $\s^{2n-1}= \re^{2n} \cap \partial \bb^{2n+1}$
which sends $(x,t) \in \bb^{2n+1}$ with $\| x \|   >0$ to
the point $(\frac {x}{\| x \| }, 0)$ in $\s^{2n-1}$. We described in detail
the basic surgery of $L^{[2n-1]}_m$ in the covering model 
because, as we mentioned before,
this description will be used in Section \ref{even-reduction-section}.

\begin{proposition}
\label{identity-odd}
The identity map of $L_m^{[2n-1]}, n\geq 1,$ extends
to a map $\alpha : M \lo L_m^{[2n]}$  from  the basic modification $M$ of $L_m^{[2n+1]}$ 
to $L_m^{[2n]}$ so that 
 $\alpha$ restricted to
$S^0\times S^{2n-1}\subset  S^1 \times S^{2n-1}=
 \partial E_M$ is null-homotopic where $S^0$ is a singleton in $S^1$.
\end{proposition}
{\bf Proof.}
Recall that $\partial E=S^1\times S^{2n-1}$ and
 $r^1$ restricted
to $S^1 \times S^0 $ and $S^0 \times S^{2n-1}$  and followed by
the inclusion of $L^{[2n-1]}_m$ into $L^{[2n]}_m$
represent a generator  of the fundamental
group of $L^{[2n]}_m$ and 
a null-homotopic map to $L^{[2n]}_m$
 respectively. 
 Then,
 By Proposition \ref{extension-odd},
the map $r^1$ restricted to $\partial E=\partial E_M$
extends over $E_M$ as a map to $L^{[2n]}_m$
and this extension together with
$r^1$ restricted to 
$L^{[2n+1]}_m \setminus (E \setminus \partial E)$ provides 
the map required  in the proposition.
$\black$
\\\\
{\bf Proof of Proposition \ref{odd-reduction-proposition}.}
Consider the basic modification $M$ of $L^{[2n+1]}_m$.
 By Theorem \ref{dranishnikov-extension-theorem},
$\edim X \leq M(\z_m, 2n-1)$.  Recall that
$\partial E=
\partial E_M \subset E_M = S^1 \times M(\z_m,2n-1)$.
Then
 $f$ restricted to $f^{-1}(\partial E)$ and followed by
the projection of $E_M=S^1 \times M(\z_m, 2n-1)$ to
$M(\z_m, 2n-1)$ extends over 
$f^{-1}(E)$ as   a map to $M(\z_m,2n-1)$ and hence 
$f$ restricted to $f^{-1}(\partial E)$ extends over 
$f^{-1}(E)$ as 
 a map to $E_M$.  The last extension together with $f$ provides
a map $f _M : X \lo M$ which coincides with $f$ on $f^{-1}(L^{[2n-1]}_m)$.
By Proposition \ref{identity-even}, take
a map  $\alpha : M \lo L^{[2n]}_m$ which extends the identity map of
$L^{[2n-1]}_m$. Then $f'=\alpha \circ f_M : X \lo L^{[2n]}$
is the map required in the proposition.
$\black$

\end{section}

\begin{section}{Pushing maps off the $(2n+2)$-skeleton of  $L_m$ }
\label{even-reduction-section}
In this section we will prove
\begin{proposition}
\label{even-reduction-proposition}
Let $X$ be a metrizable space with $\dim_{\z_m} X \leq 2n, n\geq 1$,
and let $f : X \lo  L^{[2n+2]}_m$
be  a map. Then there is a  map 
$f' : X \lo  L^{[2n+1]}_m$ which coincides with $f$ on $f^{-1}(L^{[2n]}_m)$.
\end{proposition}
The proof of Proposition \ref{even-reduction-proposition} 
is based on a modification
of $L^{[2n+2]}_m$. This modification is defined for $n\geq 1$ and
will be referred to as the basic modification of $L^{[2n+2]}_m$.
Consider the ball model 
of $L^{[2n+2]}_m$ with the projection 
$p : \bb^{2n+2} \lo L^{[2n+2]}_m$ where $\bb^{2n+2}$ is the unit
ball in $\re^{2n+2}=\re^{2n} \oplus \re^2_\perp$.
Denote by $\s^{2n+1}$,  $\s^1_\perp$  and $\bb^2_\perp$
the unit sphere, the unit circle and the unit ball
in $\re^{2n+2}$ and  $\re^2_\perp$ respectively.
Also denote $L^{[2n+1]}_m =p(\s^{2n+1})$ and
$(L^{[2]}_m)_\perp=p(\bb^2_\perp)$, and 
let $L^{[2n]}_m  \subset L^{[2n+1]}_m$ be the $(2n)$-skeleton 
of $L^{[2n+1]}_m$ 
constructed as described in Section \ref{preliminaries}.

Consider the construction of the basic modification
of $L^{[2n+1]}_m$ as described in Section \ref{odd-reduction-section}.
Extend the neighborhood  $E_\s^1$ of $\s^1_\perp$ in $\s^{2n+1}$
 to a neighborhood
$E_\bb^2$ of $\bb^2_\perp$ in $\bb^{2n+2}$ and extend the $(2n)$-ball 
$SO(2n)$-bundle structure
of $E_\s^1$ over $\s^1_\perp$ to a $(2n)$-ball $SO(2n)$-bundle
structure of $E_\bb^2$ over $\bb^2_\perp$ so that
the transformations of $\z_m$ on $\re^{2n+2}$ induce
bundle maps of $E_\bb^2$. Then the neighborhood  $E^2=p(E_\bb^2)$ of
$(L^{[2]}_m)_\perp=p(\bb^2_\perp)$ in $L^{[2n+2]}_m$
 will have
the induced $(2n)$-ball $SO(2n)$-bundle structure over
$(L^{[2]}_m)_\perp$. The retraction 
$r_\s^1 : \s^{2n+1} \setminus \s^1 \lo \s^{2n-1}$ can be extended
to a  retraction 
$r_\bb^2 : \bb^{2n+2} \setminus \bb^2_\perp \lo \s^{2n-1}$
which induces the retraction 
$r^2 : L^{[2n+2]}\setminus (L^{[2]}_m)_\perp \lo L^{[2n-1]}_m$.

 Note that 
$(L^{[0]}_m)_\perp=L^{[2n]}_m \cap (L^{[2]}_m)_\perp$ is
a singleton lying in $(L^{[1]}_m)_\perp =p(\s^1_\perp)=S^1$.
Also note that the pair $((L^{[2]}_m)_\perp,(L^{[1]}_m)_\perp)$
is homeomorphic to the pair 
$(L^{[2]}_m,L^{[1]}_m)$.
In order to simplify  the notation, 
from now we will write 
$L^{[0]}_m, L^{[1]}_m, L^{[2]}_m$ instead of 
$(L^{[0]}_m)_\perp, (L^{[1]}_m)_\perp, (L^{[2]}_m)_\perp$
 keeping in mind
that any skeleton whose dimension does not depend on $n$ 
should be interpreted as having the subscript $\perp$.

Let $E^0=E^2\cap L^{[2n]}_m $ be
the fiber of the bundle $E^2$ over the point
$L^{[0]}_m$. Denote by $\partial E^0$ the boundary of the ball $E^0$
and 
 by $\partial E^2$ the induced
$S^{2n-1}$-bundle formed by the boundaries of
the fibers of $E^2$ which are $(2n)$-balls. Consider
the retraction $r^2$ as a map to $L^{[2n-1]}_m$ 
followed by the inclusion into $L^{[2n]}_m$.
Note that $r^2$ extends $r^1$ and recall
that  $r^1$  is a deformation retraction on
 $L^{[2n]}_m \setminus L^{[1]}_m$. Then
$r^2$ can be  homotoped
into a map  $r^2_*: L^{[2n+2]}_m \setminus L^{[2]}_m\lo L^{[2n]}_m$ 
which does not move the points of 
$L^{[2n]}_m \setminus L^{[2]}_m=L^{[2n]}_m \setminus L^{[1]}_m$.
Thus we can define the map 
$r^2_0 : (L^{[2n+2]}_m \setminus (E^2 \setminus \partial E^2))\cup E^0
\lo L^{[2n]}_m$ which coincides with $r^2_*$ on
$L^{[2n+2]}_m \setminus (E^2 \setminus \partial E^2)$ and
does not move the points of $E^0$. Take a neighborhood $U$
of $E^0$ in $E$ and extend $r^2_0$ to a map
$r^2_U :  
(L^{[2n+2]}_m \setminus (E^2 \setminus \partial E^2))\cup U
\lo L^{[2n]}_m$. Consider separately the quotient space
$E=E^2/E^0$ and
consider 
$\partial E =\partial E^2/\partial E^0$ as 
as a subspace of $E$.
By Proposition \ref{embedding} embed the space $E$ into
$L^{[2n+2]}_m$ so that $E \subset E^2$,
$E\cap L^{[2n]}=L^{[0]}_m=$ the singleton $E^0$ in $E$ and
$E^2 \setminus E\subset U$. Thus we have that $r^2_U$ is defined 
on $L^{[2n+2]}_m \setminus (E \setminus \partial E)$,
$L^{[2n]}_m \subset L^{[2n+2]}_m \setminus (E \setminus \partial E)$
and hence $r^2_U$ acts on 
$L^{[2n+2]}_m \setminus (E \setminus \partial E)$ as a retraction to
$L^{[2n]}_m$.

The basic modification of $L^{[2n+2]}_m, n\geq 1$, 
is defined as follows.
 By Proposition \ref{trivial-bundle}
represent 
$\partial E=\partial E^2/\partial E^0$ as 
$$ \partial E=
\frac{L^{[2]}_m\times S^{2n-1}}{ L^{[0]}_m\times S^{2n-1}}. $$
Denote

$$E_M=
\frac{L^{[2]}_m\times M(\z_m, {2n-1})}{ L^{[0]}_m\times M(\z_m,{2n-1})},
\partial E_M=
\frac{L^{[2]}_m\times \partial M(\z_m, {2n-1})}{ L^{[0]}_m
\times \partial M(\z_m,{2n-1})}$$
and consider $\partial E_M$ as a subset of $E_M$.
By the {\bf basic surgery }of $L^{[2n+2]}_m$  we mean
 replacing $E$ with $E_M$
such that $\partial E_M$ is identified with $\partial E$
through an identification of $\partial M(\z_m, 2n-1)$
with $S^{2n-1}$.  By the {\bf basic modification} $M$ of $L^{[2n+2]}_m$
we mean
the space obtained from $L^{[2n+2]}_m$ by the basic surgery
of $L^{[2n+2]}_m$.
Clearly $L^{[2n]}_m$ remains
untouched in  $M$.
\begin{proposition}
\label{identity-even}
The identity map of $L_m^{[2n]}, n\geq 1,$ extends
to a map from  the basic modification $M$ of $L_m^{[2n+2]}$ 
to $L_m^{[2n+1]}$.
\end{proposition}
{\bf Proof.} By Proposition \ref{extension-even}
the map $r^2_U$ restricted to $\partial E=\partial E_M$
extends over $E_M$ and this extension together with
$r^2_U$ restricted to 
$L^{[2n+2]}_m \setminus (E \setminus \partial E)$ provides 
the map required  in the proposition.
$\black$
\\\\
{\bf Proof of Proposition \ref{even-reduction-proposition}.}
Consider the basic modification $M$ of $L^{[2n+2]}_m$.
 By Theorem \ref{dranishnikov-extension-theorem},
$\edim X \leq \Sigma M(\z_m, 2n-1)$.  Recall that
$\partial E=
\partial E_M \subset E_M =
\frac{L^{[2]}_m \times M(\z_m,2n-1)}{L^{[0]}_m \times M(\z_m, 2n-1)}$.
Then, by  Proposition \ref{suspension},
 $f$ restricted to $f^{-1}(\partial E)$ extends over 
$f^{-1}(E)$ as a map to $E_M$ and this extension together with $f$ provides
a map $f _M : X \lo M$ which coincides with $f$ on $f^{-1}(L^{[2n]}_m)$.
By Proposition \ref{identity-even}, take
a map  $\alpha : M \lo L^{[2n+1]}_m$ which extends the identity map of
$L^{[2n]}_m$. Then $f'=\alpha \circ f_M : X \lo L^{[2n+1]}$
is the map required in the proposition.
$\black$

\end{section}

\begin{section}{Pushing maps off the $3$-skeleton of $L_m$}
\label{reduction-in-dim-3}
The goal of this section is to prove Theorem \ref{main-theorem}.
Clearly  Propositions \ref{odd-reduction-proposition}
and \ref{even-reduction-proposition}  imply
\begin{theorem}
\label{reduction-to-dim-3}
Let $X$ be a finite dimensional metrizable space with $\dim_{\z_m} X \leq 2$
and let $f : X \lo  L^{[n]}_m, n\geq 3,$ 
be  a map. Then there is a  map 
$f' : X \lo  L^{[3]}_m$ which coincides with $f$ on $f^{-1}(L^{[2]}_m)$.
\end{theorem}
An easy corollary of Theorem \ref{reduction-to-dim-3} is 
\begin{corollary}
\label{corollary-dim-3}
Let $X$ be a finite dimensional metrizable space with
$\dim_{\z_m} X \leq 1$ and $f_F : F \lo L^{[2]}_m$ 
a map from a closed subset $F$ of $X$. Then $f_F$ extends
to a map $ f : X \lo L^{[3]}_m$.
\end{corollary}
{\bf Proof.}
Since $L_m=K(\z_m,1)$ 
we have
$\edim X \leq L_m$.  Then $f_F$ extends
to a map $f : X \lo L_m$. Since $X$ is finite dimensional 
 we can assume that there is $n$ such that $f( X )\subset  L^{[n]}_m$. 
 Then, by Theorem \ref{reduction-to-dim-3}, 
 one can replace $f$ by a map to $L^{[3]}_m$ which coincides  
 with $f_F$ on $F$  and
 the corollary follows.
 $\black$
 \\\\
 Thus  the only missing part of proving Theorem \ref{main-theorem}
 is to push maps from  $L^{[3]}_m$ to $L^{[2]}_m$. We will do this in
 two steps. The first one is

\begin{proposition}
\label{reduction-in-dim-3-first}
Let $X$ be a finite dimensional 
metrizable space with $\dim_{\z_m} X \leq 1$ and
let $f : X \lo  L^{[3]}_m$ 
be  a map. Then there is a  map 
$f' : X \lo  L^{[2]}_m$ which coincides with $f$ on $f^{-1}(L^{[1]}_m)$.
\end{proposition}
{\bf Proof.} 
Consider the basic modification $M$ of
$L^{[3]}_m$. 
Recall that $M$ is obtained from $L^{[3]}_m$ 
by the basic surgery which replaces
 $E=S^1 \times B \subset L^{[3]}_m$
with $E_M =S^1 \times M(\z_m,1)=S^1 \times L^{[2]}_m$
by identifying 
$\partial E=S^1 \times \partial B =S^1 \times S^1$
with $\partial E_M=S^1 \times  L^{[1]}_m=
S^1\times S^1$.
Also recall that $L^{[1]}_m$ remains untouched  in $M$ 
and does not meet $E_M$.

Enlarge $M$ to the space $M^+$ by enlarging
$E_M = S^1 \times L^{[2]}_m$ to $E^+_M = S^1 \times L^{[3]}_m$.
Apply again  Corollary \ref{corollary-dim-3} to 
 the projection  of $S^1 \times L^{[3]}_m$ to $L^{[3]}_m$
to extend the map
 $f$ restricted to $f^{-1}(\partial E)$  over
  $f^{-1}(E)$ as  a map to $S^1 \times L^{[3]}_m$ and this way
 to get a map $f^+ : X \lo M^+$ which differs from $f$
 only on $f^{-1}(E)$.

The space $M^{++} $ is obtained from $M^+$ by
replacing $E^+_M=S^1 \times L^{[3]}_m$ with 
$E^{++}_M=S^1\times M$ by identifying
$L^{[1]}_m $ in $L^{[3]}_m$ with $L^{[1]}_m $ in  $M$. 
Thus $M^{++}$ differs from $M^+$ on the set 
$E^{++}_M =S^1 \times S^1 \times L^{[2]}_m\subset M^{++}$.

By Proposition \ref{identity-odd} the identity map 
of $L^{[1]}_m$ extends to  a map 
$\alpha : M \lo L^{[2]}_m$ which is 
null-homotopic on 
$S^0 \times L^{[1]}_m \subset E_M$ where   $S^0$ is a singleton in $ S^1$.
Then we get that the map
$id \times \alpha : 
E^{++}_M= S^1 \times M \lo S^1\times  L^{[2]}_m \subset E^+_M$
induces 
 $\beta : M^{++} \lo M^{+}$
so that $\beta(M^{++})\subset M \subset M^+$.
Consider the map 
$\gamma=\alpha \circ \beta : M^{++} \lo L^{[2]}_m$
and note that $\gamma$ restricted to 
$S^0 \times S^0 \times L^{[2]}_m \subset 
E^{++}_M$ is the composition of the maps
$S^0 \times S^0 \times L^{[2]}_m \lo S^0 \times L^{[2]}_m \lo L^{[2]}_m$
each of them acting as   $\alpha $ restricted to 
$S^0 \times L^{[2]}_m\subset E_M$.
 Hence, by Proposition \ref{composition}, 
$\gamma$ restricted to 
$S^0 \times S^0 \times L^{[2]}_m $
is null-homotopic. Then, by Proposition \ref{suspension},
$f^+$ restricted to $(f^+)^{-1}(\partial E^{++}_M))$ and followed
by $\gamma$ for 
 $\partial E^{++}_M =S^1\times S^1 \times L^{[1]}_m
\subset E^{++}_M$  extends over $(f^+)^{-1}( E^{++}_M))$
and this extension provides
 a map $f' : X \lo L^{[2]}_m$ that coincides
with $f$ on $f^{-1}(L^{[1]}_m)$. The proposition is proved.
$\black$.
\begin{proposition}
\label{reduction-in-dim-3-second}
Let $X$ be a finite dimensional 
metrizable space with $\dim_{\z_m} X \leq 1$ and
let $f : X \lo  L^{[3]}_m$ 
be  a map. Then there is a  map 
$f' : X \lo  L^{[2]}_m$ which coincides with $f$ on $f^{-1}(L^{[2]}_m)$.
\end{proposition}
{\bf Proof.} Consider the ball model of $L^{[3]}_n$ and
let $p : \bb^3 \lo L^{[2]}_m$.
In the decomposition $\re^3=\re^2 \oplus \re_\perp$ we will
refer to $\re^2$ and $\re_\perp $
as the $xy$-coordinate plane and  the $z$-axis respectively.
By a rotation of $\bb^3$ we mean an orthogonal
 rotation around the $z$-axis.
Clearly 
a rotation $\phi$ of $\bb^3$ induces  a  homeomorphism $\phi_L$ of
$L^{[3]}_m$ which  will be  called the induced rotation  
of $L^{[3]}_m$. Note that  the rotations of $\bb^3$ and $L^{[3]}_m$
commute  with the projection $p$. 
 By this we mean that 
$p \circ \phi=\phi_L \circ p$.
Also note that $L^{[2]}_m$ is invariant under
the induced rotations of $L^{[3]}_m$.

Take a disk $B$ of radius $1/3$ lying in the $xz$-coordinate
plane and centered at the point $(1/2,0,0) \in \re^3$.
Denote by $E=S^1 \times B$ the solid torus obtained by rotating
$B$ around the $z$-axis and  denote 
by $\partial E=S^1 \times \partial B=S^1 \times S^1\subset E$ 
the boundary
of $E$. Clearly $E$  can be considered as subsets of $L^{[3]}_m$.

  Let  $\I =\bb^3 \cap \re_\perp  $ be the $[-1,1]$-interval
of the $z$-axis
and $\partial \I =\{-1, 1\}\subset \re_\perp$ 
the end points of $\I$.  Consider an obvious   retraction
$r_\bb : \bb^3 \setminus (E \setminus \partial E ) 
\lo \I\cup \partial \bb^3=\I \cup \s^2$ such that $r_\bb$ commutes with 
the rotations of $\bb^3$  and
consider the map $\gamma: \I \cup \partial \bb^3 \lo L^{[2]}_m$
such that 
$\gamma$  coincides on $\s^2$  with $p$
and $\gamma$ sends $\I$ to the point $p(\partial \I)$.
Note that $r_\bb$ followed by $\gamma$ induces the retraction
$r : L^{[3]}_m \setminus (E \setminus \partial E) \lo L^{[2]}_m$ such that
$r$ commutes with the rotations of $L^{[3]}_m$.

Consider the surgery of $L^{[3]}_m$  which replaces $E=S^1\times B$
with $E_M = S^1 \times L^{[2]}_m$ by identifying the boundary 
$\partial B=S^1$
 of $B$
with the $1$-skeleton $L^{[1]}_m=S^1$ of $L^{[2]}_m$. 
Denote by $M$ the space obtained from $L^{[3]}_m$ by this surgery.
Clearly the $2$-skeleton $L^{[2]}_m$ of $L^{[3]}_m$ remains untouched 
in $M$.
Note that   this surgery and this modification are
 different 
from the basic surgery and the basic modification of $L^{[3]}_m$
considered before.

Observe that any map from $L^{[1]}_m$ to $L^{[2]}_m$
extends over  $L^{[2]}_m$.
Fix a singleton $S^0$
in $S^1$ and extend $r$ restricted to 
$S^0 \times L^{[1]}_m$  to a map
$r^0: S^0 \times L^{[2]}_m \lo L^{[2]}_m$.
Then, since $r$ commutes with the rotations of $L^{[3]}_m$,
$r^0$ can be extended by the rotations of $L^{[3]}_m$
to a map $r^1 :E_M= S^1 \times L^{[2]}_m \lo L^{[2]}_m$
so that $r^1$ extends $r$ restricted to $\partial E$.
Thus $r^1$ together with $r$  provide a retraction $\beta : M \lo L^{[2]}_m$
which extends the identity map of $L^{[2]}_m$.

Now denote  $F=f^{-1}(\partial E)$ and consider the restriction
$f|F : F \lo \partial E =S^1 \times S^1=S^1 \times L^{[1]}_m$.
Then, by Corollary \ref{corollary-dim-3} and
Proposition \ref{reduction-in-dim-3-first},
 $f|F$ followed by the projection to $L^{[1]}_m$ extends
 over $f^{-1}(E)$ as a map to $L^{[2]}_m$ and this extension provides
 an extension  $f_E : f^{-1}(E) \lo E=S^1 \times L^{[2]}_m$
 of $f|F $ over $f^{-1}(E)$. Thus we get 
 a map $f_M : X \lo M$  which
 coincides with $f$  on $f^{-1}(L^{[2]}_m)$. Set 
 $f'=\beta\circ f_M : X \lo L^{[2]}_m$ and the proposition follows.
 $\black$
 \\
\\
{\bf Proof of Theorem \ref{main-theorem}. }Theorem \ref{main-theorem}
 follows
from Corollary \ref{corollary-dim-3} and Proposition 
\ref{reduction-in-dim-3-second}. $\black$

\end{section}

Jerzy Dydak \\                 
 Department of Mathematics\\
 University of Tennessee\\
 Knoxville, TN 37996-1300 \\                 
 dydak@math.utk.edu 
\\\\ 
Michael Levin\\
Department of Mathematics\\
Ben Gurion University of the Negev\\
P.O.B. 653\\
Be'er Sheva 84105, ISRAEL  \\
 mlevine@math.bgu.ac.il\\\\

\begin{thebibliography}{99}



\bibitem{DrExt} A.N. Dranishnikov {\em An extension of mappings into CW-complexes}, 
Mat. Sb. 182 (1991), 1300-1310; English transl., Math. USSR Sb. 74 (1993), 47-56.


\bibitem{DrArxiv} Dranishnikov, A. N. {\em Cohomological dimension theory of compact 
metric spaces,} Topology
Atlas invited contribution, http://at.yorku.ca/topology.taic.html 
(see also  arXiv:math/0501523).

\bibitem{Dy}
 Dydak, Jerzy {\em Cohomological dimension and metrizable spaces. II.} 
Trans. Amer. Math. Soc. 348 (1996), no. 4, 1647--1661.
\bibitem{hatcher} 

Hatcher, Allen {\em Algebraic topology}. Cambridge University Press, Cambridge,
 2002. xii+544 pp. ISBN: 0-521-79160-X; 0-521-79540-0
\bibitem{first-projective}
Dydak, Jerzy; Levin, Michael {\em Extensions of maps to the projective
  plane}.  Algebr. Geom. Topol.,
  5(2005), 1711-1718.

\bibitem{dydak-levin-2}
Dydak, Jerzy; Levin, Michael {\em  Maps to the projective plane}.  
Algebr. Geom. Topol.  9  (2009),  no. 1, 549-568. 
 
\bibitem{Ku}
Kuzminov, V. I. {\em Homological dimension theory.} Russian Math Surveys 23 (5)  (1968), 
1-45.










\end{thebibliography}
\end{document}